\documentclass[11pt]{article}
\usepackage[margin=1in]{geometry}
\usepackage{amsmath,amssymb,amsthm,tikz}
\usepackage{amssymb,authblk}
\usepackage[hidelinks]{hyperref}
\newtheorem{theorem}{Theorem}[section]

\newcommand{\R}{\mathbb R}
\newcommand{\E}{\mathbb E}
\newcommand{\Prob}{\mathbb P}
\newcommand{\eps}{\varepsilon}
\title{On the measure of Neighborhoods of Kakeya and Truncated Nikodym Sets}
\author{Iqra Altaf}
\date{}
\begin{document}
\maketitle
\begin{abstract}
We show that there exist Kakeya sets in $\R^d$ whose $\delta$-neighborhoods have measure $\lesssim |\log\delta|^{-(d-1)}$.
We also show that the measure of a $\delta$-neighborhood of an $\eps$-truncated Nikodym set in $\R^2$ $\gtrsim \eps^{-1}|\log\delta|^{-1}$, for $\delta$ small enough. The sufficiently small scale depends on the Nikodym set and on $\eps$.
\end{abstract}

\section{Introduction}
A Kakeya set in $\R^d$ is a null set which contains a unit line segment in every direction \cite{Bes}. In the plane it has long been known that Kakeya sets have Hausdorff dimension $2$. The Kakeya set conjecture states that every Kakeya set in $\R^d$ has Hausdorff dimension $d$. In a remarkable breakthrough, Wang and Zahl proved the Kakeya conjecture in $\R^3$ \cite{WZ3}. Their proof distinguishes between sticky Kakeya sets and non-sticky Kakeya sets \cite{WZ1,WZ2,WZ3}. In \cite{LW}, Lai and Wong constructed a non-sticky Kakeya set of Lebesgue measure zero in the plane. In higher dimensions, they constructed sticky and non-sticky Kakeya sets which are not Cartesian products of a planar Kakeya set with $\R^{d-2}$.

For a Kakeya set in $\R^2$, a standard argument shows that its $\delta$-neighborhood must have measure $\gtrsim |\log\delta|^{-1}$ \cite{Cor}. In \cite{Keich}, Keich showed that this bound is sharp by explicitly constructing a Kakeya set whose $\delta$-neighborhood has measure $\lesssim |\log\delta|^{-1}$. Therefore, by taking a Cartesian product of $d/2$ such Kakeya sets (for even $d$), we get a Kakeya set in $\R^d$ whose $\delta$-neighborhood has measure $\lesssim_d |\log\delta|^{-d/2}$. Taking a Cartesian product of this Kakeya set with a unit interval, we get a Kakeya set in $\R^{d+1}$ whose $\delta$-neighborhood also has measure $\lesssim_d |\log\delta|^{-d/2}$. In this paper, we construct Kakeya sets that have $\delta$-neighborhoods of measure $\lesssim_d |\log\delta|^{-(d-1)}$. We do this by building on the recent work of Lai and Wong \cite{LW}. Fern\'andez-Delgado and de la Salle obtained the same bound by a different construction \cite{FdS}.

Closely related to Kakeya sets are Nikodym sets. A Nikodym set $N$ in $\R^d$ is a null set such that for every $y\in\R^d$ there is a line $\ell_y$ through $y$ for which the punctured line $\ell_y\setminus\{y\}$ is contained in $N$. The original planar construction is due to Nikodym \cite{Nik}.

Falconer \cite{Fal} used sets with prescribed projections
to obtain a higher-dimensional version of Nikodym's construction.
He constructed a null set $F\subset\R^d$ such that through every point
$x\notin F$ there is a hyperplane $H_x$ with
$H_x\setminus\{x\}\subset F$. In this paper we return to line segments
and ask how small a neighborhood of their union can be after we remove
the part of each segment nearest its base point.

\paragraph{Notation.} For $d\in\mathbb N$, a set $E\subset\R^d$ and $\delta>0$, $|E|$ denotes the $d$-dimensional Lebesgue measure of $E$ and $E(\delta)$ denotes the $\delta$-neighborhood of $E$, i.e.,
\[
 E(\delta)=\{x\in\R^d:\operatorname{dist}(x,E)<\delta\}.
\]
Our first theorem is about Kakeya sets.
\begin{theorem}\label{thm:k}
For every integer $d\ge2$, there exists a compact Kakeya set $K\subset\R^d$ such that
\[
 |K(\delta)|\lesssim_d (\log(1/\delta))^{-(d-1)}
\]
for all sufficiently small $\delta>0$.
\end{theorem}

We actually prove more than Theorem~\ref{thm:k}. We will construct Kakeya sets from random Cantor sets. We show that the $\delta$-neighborhoods of almost all such Kakeya sets in $\R^d$ have measure $\lesssim_d (\log(1/\delta))^{-(d-1)}$.

For the Nikodym statement, let $z:[0,1]^2\to\R$ be a Lebesgue-measurable function. Let $N\subset\R^2$ be a null set such that
\begin{equation}\label{inc}
N=\{ p+t(z(p),1) :\text{for every }p\in [0,1]^2\text{ and every }0<t\le1\}.
\end{equation}
Clearly, such a set $N$ is dense along $[0,1]^2$, in the sense that $[0,1]^2\subset\overline N$. Therefore, in order to meaningfully ask the question about the measure of the $\delta$-neighborhood of $N$, we consider a truncated version of a Nikodym set. Define the $\eps$-truncation of $N$ as
\begin{equation}\label{Trunc}
 N_\eps=\{p+t(z(p),1):p\in [0,1]^2,\ \eps\le t\le1\}.
\end{equation}

\begin{theorem}\label{thm:n}
Let $N$ be a Nikodym set as in \eqref{inc}. There is an absolute constant $c>0$ such that, for every $0<\eps\le1/2$,
\[
 \liminf_{\delta\downarrow0}\log(1/\delta)\,|N_\eps(\delta)|\ge c/\eps .
\]
Consequently, after decreasing the absolute constant if necessary,
\[
 |N_\eps(\delta)|\ge\frac{c}{\eps\log(1/\delta)}
 \qquad(0<\delta<\delta_0(N,z,\eps)).
\]
\end{theorem}
We allow the threshold $\delta_0$ to depend on the fixed Nikodym set $N$, the choice of slopes, and $\eps$.

\section{The Kakeya construction}\label{sec:cantor}
We first describe a well-known construction of Kakeya sets via purely unrectifiable sets. Let $E\subset[0,1]^2$ be a compact, purely unrectifiable set such that $\mathcal H^1(E)<\infty$ and its projection onto the $x$-axis is $[0,1]$. Let 
\[
 K=\{(tx+y,t):(x,y)\in F,\, 0\le t\le1\}.
\]
For any $x\in[0,1]$, the projection assumption gives at least one $y$ with $(x,y)\in E$. The corresponding line segment $(tx+y,t)=(y,0)+t(x,1)$ has direction $(x,1)$ and length $\sqrt{1+x^2}\geq 1$. These directions fill an arc of the unit circle, so finitely many rotated copies of $K$ have a line segment in every direction. Let, $P_t : \mathbb{R}^2 \rightarrow \mathbb{R}$ such that
$$P_t(x,y)=tx+y.$$
By the Besicovitch-Federer projection theorem, almost every projection of a purely unrectifiable set with finite $\mathcal H^1$ has length zero. The map $P_t$ gives the scalar projection in the direction $(t,1)$ multiplied by $\sqrt{1+t^2}$. It follows that
$|P_t(E)|=0$ for almost every $t\in[0,1]$. Fubini's theorem gives $|K|=0$. The finite union of
rotated copies of $K$ is therefore a compact planar Kakeya set of measure zero. \\
For $d \geq 2$, let $E_1,\dots,E_{d-1}$ be purely unrectifiable sets as above. Let,
\begin{equation}\label{Kakeyad}
    K=\{((y_1,\ldots,y_{d-1})+t(x_1,\ldots,x_{d-1}),t):(x_j,y_j)\in E_j,\, 0\le t\le1\}.
\end{equation}
We therefore get a segment in direction $(x_1,\ldots,x_{d-1},1)$ for every
such $x_1,\ldots,x_{d-1} \in [0,1]^{d-1}$. These directions contain an open subset of the unit sphere and finitely many rotations cover all directions. Let $(z_1,\cdots,z_d)$ denote a point in $\mathbb{R}^d$. Since the pairs $(x_j,y_j)$ are chosen independently, the intersection of the set $K$ with the plane $z_d=t$ is $P_t(E_1)\times\cdots\times P_t(E_{d-1}).$ Thus Besicovitch-Federer theorem and Fubini's theorem give $|K|=0$ as before. To get the desired estimate for $|K(\delta)|$, we will choose random planar Cantor sets with small projection lengths. 

For the construction of Cantor sets, we start with the square $E_0=[0,1]^2$. Then subdivide each square into a $4\times4$ grid. In each of its four vertical columns, we retain one of the four squares,
chosen uniformly at random. After $n$ many steps, we have $4^{n}$ many squares of side length $4^{-n}$. Let their union be $E^n$. (See Figure 1) The set $E^n$ has full $x$-projection. Define
\begin{equation*} 
    E=\bigcap_{n=0}^{\infty}E^n.
\end{equation*}
Let $\mathcal{E}$ be the collection of all such sets $E$.\par 
We compare this construction with the random Cantor sets of Peres and
Solomyak \cite[Theorem~2.2]{PS}. They divide the unit square into four equal quarters. In each quarter they independently choose one of its four subsquares of side $1/4$, each with the same probability. (See Figure 2) For $0\le\theta\le\pi$, let $\operatorname{proj}_\theta$ denote orthogonal
projection onto the line through the origin in the $\theta$-direction. The Favard length of a planar set $A$ is
\[
 \operatorname{Fav}(A)=\int_0^\pi |\operatorname{proj}_\theta A|\,d\theta.
\]
Let $F^n$ be the union of the squares kept after $n$ steps. Peres and Solomyak show in \cite[(6.3)]{PS} that for every direction $\theta$, the expected length
of the orthogonal projection of $E^n$ satisfies
\[
 \mathbb E\bigl[|\operatorname{proj}_\theta F^n|\bigr]\lesssim 1/n,
\]
where the constant is independent of $\theta$ and $n$. Integrating over $\theta$ gives their Favard length estimate $\mathbb E[\operatorname{Fav}(F^n)]\lesssim 1/n$. Their proof applies to our
choice of one square per column for the directions $0\le t\le1$ and gives
\begin{equation}\label{mean}
 \mathbb E[|P_t(E^n)|]\lesssim 1/n,  \qquad \text{ for all, }0\le t\le1.
\end{equation}
Since the exactly same argument applies, we do not reproduce the proof here.

\begin{figure}
\centering
\begin{tikzpicture}[x=3.7cm,y=3.7cm]
\begin{scope}[shift={(0.0,0)}]
\draw[black!75,line width=.45pt] (0.25, 0.0) -- (0.25, 1.0);
\draw[black!35,line width=.2pt] (0.0, 0.25) -- (1.0, 0.25);
\draw[black!75,line width=.45pt] (0.5, 0.0) -- (0.5, 1.0);
\draw[black!35,line width=.2pt] (0.0, 0.5) -- (1.0, 0.5);
\draw[black!75,line width=.45pt] (0.75, 0.0) -- (0.75, 1.0);
\draw[black!35,line width=.2pt] (0.0, 0.75) -- (1.0, 0.75);
\fill (0.0, 0.5) rectangle (0.25, 0.75);
\fill (0.25, 0.25) rectangle (0.5, 0.5);
\fill (0.5, 0.25) rectangle (0.75, 0.5);
\fill (0.75, 0.5) rectangle (1.0, 0.75);
\draw (0,0) rectangle (1,1);
\node[below] at (.5,-.04) {Step 1};
\end{scope}
\begin{scope}[shift={(1.32,0)}]
\draw[black!75,line width=.45pt] (0.0625, 0.5) -- (0.0625, 0.75);
\draw[black!35,line width=.2pt] (0.0, 0.5625) -- (0.25, 0.5625);
\draw[black!75,line width=.45pt] (0.125, 0.5) -- (0.125, 0.75);
\draw[black!35,line width=.2pt] (0.0, 0.625) -- (0.25, 0.625);
\draw[black!75,line width=.45pt] (0.1875, 0.5) -- (0.1875, 0.75);
\draw[black!35,line width=.2pt] (0.0, 0.6875) -- (0.25, 0.6875);
\draw[black!55,thin] (0.0, 0.5) rectangle (0.25, 0.75);
\draw[black!75,line width=.45pt] (0.3125, 0.25) -- (0.3125, 0.5);
\draw[black!35,line width=.2pt] (0.25, 0.3125) -- (0.5, 0.3125);
\draw[black!75,line width=.45pt] (0.375, 0.25) -- (0.375, 0.5);
\draw[black!35,line width=.2pt] (0.25, 0.375) -- (0.5, 0.375);
\draw[black!75,line width=.45pt] (0.4375, 0.25) -- (0.4375, 0.5);
\draw[black!35,line width=.2pt] (0.25, 0.4375) -- (0.5, 0.4375);
\draw[black!55,thin] (0.25, 0.25) rectangle (0.5, 0.5);
\draw[black!75,line width=.45pt] (0.5625, 0.25) -- (0.5625, 0.5);
\draw[black!35,line width=.2pt] (0.5, 0.3125) -- (0.75, 0.3125);
\draw[black!75,line width=.45pt] (0.625, 0.25) -- (0.625, 0.5);
\draw[black!35,line width=.2pt] (0.5, 0.375) -- (0.75, 0.375);
\draw[black!75,line width=.45pt] (0.6875, 0.25) -- (0.6875, 0.5);
\draw[black!35,line width=.2pt] (0.5, 0.4375) -- (0.75, 0.4375);
\draw[black!55,thin] (0.5, 0.25) rectangle (0.75, 0.5);
\draw[black!75,line width=.45pt] (0.8125, 0.5) -- (0.8125, 0.75);
\draw[black!35,line width=.2pt] (0.75, 0.5625) -- (1.0, 0.5625);
\draw[black!75,line width=.45pt] (0.875, 0.5) -- (0.875, 0.75);
\draw[black!35,line width=.2pt] (0.75, 0.625) -- (1.0, 0.625);
\draw[black!75,line width=.45pt] (0.9375, 0.5) -- (0.9375, 0.75);
\draw[black!35,line width=.2pt] (0.75, 0.6875) -- (1.0, 0.6875);
\draw[black!55,thin] (0.75, 0.5) rectangle (1.0, 0.75);
\fill (0.0, 0.5625) rectangle (0.0625, 0.625);
\fill (0.0625, 0.5) rectangle (0.125, 0.5625);
\fill (0.125, 0.625) rectangle (0.1875, 0.6875);
\fill (0.1875, 0.5625) rectangle (0.25, 0.625);
\fill (0.25, 0.375) rectangle (0.3125, 0.4375);
\fill (0.3125, 0.375) rectangle (0.375, 0.4375);
\fill (0.375, 0.4375) rectangle (0.4375, 0.5);
\fill (0.4375, 0.3125) rectangle (0.5, 0.375);
\fill (0.5, 0.25) rectangle (0.5625, 0.3125);
\fill (0.5625, 0.25) rectangle (0.625, 0.3125);
\fill (0.625, 0.375) rectangle (0.6875, 0.4375);
\fill (0.6875, 0.25) rectangle (0.75, 0.3125);
\fill (0.75, 0.5) rectangle (0.8125, 0.5625);
\fill (0.8125, 0.625) rectangle (0.875, 0.6875);
\fill (0.875, 0.5) rectangle (0.9375, 0.5625);
\fill (0.9375, 0.625) rectangle (1.0, 0.6875);
\draw (0,0) rectangle (1,1);
\node[below] at (.5,-.04) {Step 2};
\end{scope}

\end{tikzpicture}
\caption{Our Cantor set. Keep one square in each of its four columns.}
\medskip
\centering
\begin{tikzpicture}[x=3.7cm,y=3.7cm]
\begin{scope}[shift={(0.0,0)}]
\draw[black!35,line width=.2pt] (0.25, 0.0) -- (0.25, 1.0);
\draw[black!35,line width=.2pt] (0.0, 0.25) -- (1.0, 0.25);
\draw[black!75,line width=.45pt] (0.5, 0.0) -- (0.5, 1.0);
\draw[black!75,line width=.45pt] (0.0, 0.5) -- (1.0, 0.5);
\draw[black!35,line width=.2pt] (0.75, 0.0) -- (0.75, 1.0);
\draw[black!35,line width=.2pt] (0.0, 0.75) -- (1.0, 0.75);
\fill (0.25, 0.0) rectangle (0.5, 0.25);
\fill (0.0, 0.5) rectangle (0.25, 0.75);
\fill (0.5, 0.25) rectangle (0.75, 0.5);
\fill (0.5, 0.75) rectangle (0.75, 1.0);
\draw (0,0) rectangle (1,1);
\node[below] at (.5,-.04) {Step 1};
\end{scope}
\begin{scope}[shift={(1.32,0)}]
\draw[black!35,line width=.2pt] (0.3125, 0.0) -- (0.3125, 0.25);
\draw[black!35,line width=.2pt] (0.25, 0.0625) -- (0.5, 0.0625);
\draw[black!75,line width=.45pt] (0.375, 0.0) -- (0.375, 0.25);
\draw[black!75,line width=.45pt] (0.25, 0.125) -- (0.5, 0.125);
\draw[black!35,line width=.2pt] (0.4375, 0.0) -- (0.4375, 0.25);
\draw[black!35,line width=.2pt] (0.25, 0.1875) -- (0.5, 0.1875);
\draw[black!55,thin] (0.25, 0.0) rectangle (0.5, 0.25);
\draw[black!35,line width=.2pt] (0.0625, 0.5) -- (0.0625, 0.75);
\draw[black!35,line width=.2pt] (0.0, 0.5625) -- (0.25, 0.5625);
\draw[black!75,line width=.45pt] (0.125, 0.5) -- (0.125, 0.75);
\draw[black!75,line width=.45pt] (0.0, 0.625) -- (0.25, 0.625);
\draw[black!35,line width=.2pt] (0.1875, 0.5) -- (0.1875, 0.75);
\draw[black!35,line width=.2pt] (0.0, 0.6875) -- (0.25, 0.6875);
\draw[black!55,thin] (0.0, 0.5) rectangle (0.25, 0.75);
\draw[black!35,line width=.2pt] (0.5625, 0.25) -- (0.5625, 0.5);
\draw[black!35,line width=.2pt] (0.5, 0.3125) -- (0.75, 0.3125);
\draw[black!75,line width=.45pt] (0.625, 0.25) -- (0.625, 0.5);
\draw[black!75,line width=.45pt] (0.5, 0.375) -- (0.75, 0.375);
\draw[black!35,line width=.2pt] (0.6875, 0.25) -- (0.6875, 0.5);
\draw[black!35,line width=.2pt] (0.5, 0.4375) -- (0.75, 0.4375);
\draw[black!55,thin] (0.5, 0.25) rectangle (0.75, 0.5);
\draw[black!35,line width=.2pt] (0.5625, 0.75) -- (0.5625, 1.0);
\draw[black!35,line width=.2pt] (0.5, 0.8125) -- (0.75, 0.8125);
\draw[black!75,line width=.45pt] (0.625, 0.75) -- (0.625, 1.0);
\draw[black!75,line width=.45pt] (0.5, 0.875) -- (0.75, 0.875);
\draw[black!35,line width=.2pt] (0.6875, 0.75) -- (0.6875, 1.0);
\draw[black!35,line width=.2pt] (0.5, 0.9375) -- (0.75, 0.9375);
\draw[black!55,thin] (0.5, 0.75) rectangle (0.75, 1.0);
\fill (0.3125, 0.0625) rectangle (0.375, 0.125);
\fill (0.25, 0.1875) rectangle (0.3125, 0.25);
\fill (0.4375, 0.0625) rectangle (0.5, 0.125);
\fill (0.375, 0.1875) rectangle (0.4375, 0.25);
\fill (0.0625, 0.5) rectangle (0.125, 0.5625);
\fill (0.0625, 0.625) rectangle (0.125, 0.6875);
\fill (0.125, 0.5) rectangle (0.1875, 0.5625);
\fill (0.125, 0.625) rectangle (0.1875, 0.6875);
\fill (0.5, 0.3125) rectangle (0.5625, 0.375);
\fill (0.5, 0.4375) rectangle (0.5625, 0.5);
\fill (0.625, 0.25) rectangle (0.6875, 0.3125);
\fill (0.6875, 0.375) rectangle (0.75, 0.4375);
\fill (0.5, 0.8125) rectangle (0.5625, 0.875);
\fill (0.5625, 0.9375) rectangle (0.625, 1.0);
\fill (0.625, 0.75) rectangle (0.6875, 0.8125);
\fill (0.625, 0.875) rectangle (0.6875, 0.9375);
\draw (0,0) rectangle (1,1);
\node[below] at (.5,-.04) {Step 2};
\end{scope}

\end{tikzpicture}
\caption{The Peres-Solomyak Cantor set. Keep one small square in each of its four quarters.}
\end{figure}
\begin{proof}[Proof of Theorem~\ref{thm:k}]
Take independent cantor sets $E_1,\dots,E_{d-1} \in \mathcal{E}$ and define the set $K$ as in \eqref{Kakeyad}.

Write $X_{j,n}(t)=|P_t(E_j^n)|$. Take $4^{-n-1}<\delta\le4^{-n}$. Let $(z_1, \cdots z_{d-1},t)\in K(\delta)$ and $0\le t\le1$, then $|z_j-(t'x_j+y_j)| \leq \delta$ for some $(x_j,y_j)\in E_j$ and $|t'-t|\leq \delta$. Since $x_j \in [0,1]^2$, $z_j$ lies within $2\delta$ of $P_t(E_j) \subset P_t(E_j^n)$. The set $P_t(E_j^n)$ contains projections of square of length $4^{-n}$. Hence each maximal interval in $P_t(E_j^n)$ has length at least $4^{-n}$. Adding $2\delta$ at each end increases its
length by at most a factor of five. The part of $K(\delta)$ such that $z_{d} \neq [0,1]$
has volume $O(\delta)$ since $K\subset[0,2]^{d-1}\times[0,1]$. Therefore
\begin{equation}\label{slice}
 |K(\delta)|\le O(\delta)+5^q\int_0^1\prod_{j=1}^q |P_t(E_j^n)|\,dt .
\end{equation}
Let 
$$ I_{n}:= \int_0^1\prod_{j=1}^q |P_t(E_j^n)|\,dt .$$
Using \eqref{mean} and independence of the sets $E_{j}$, we get that the expectation of $I_n$ is  $\lesssim (1/n)^{d-1}$. Therefore the expected value of $|K(\delta)| \lesssim |\log \delta|^{-(d-1)}$. We want to establish that for some fixed set $K,\, |K(\delta)|$ has the required bound for all sufficiently small $\delta$.\par 
Fix $t$. Suppose we have chosen the squares up to step $n$ in each of the $d-1$ sets $E_j$ such that $I_n \leq (2C/n)^{d-1}$. We will now choose the squares from step $n+1$ through
step $2n$. For each $E_j$, the $4^{n}$ squares are  filled independently. So we are making the same random choices but at scale $4^{-n}$. Therefore in each such square the expected length of the projection after these $n$ steps is at most $\sim 4^{-n}/n$. Since there are $4^{n}$ squares we get 
\begin{equation}\label{expectation}
    \E[|P_t(E_j^{2n})| \,| E_j^{n} \text{ is given}]\le 4^n\cdot 4^{-n}\frac Cn=\frac Cn.
\end{equation} 
By independence of the sets $E_j$, we have 
$$  \E\!\left[\prod_{j=1}^q |P_t(E_j^{2n})|\, E_j^{n} \text{ is given for all }j\right]\, 
 =\prod_{j=1}^{d-1}\E[|P_t(E_j^{2n})| \,| E_j^{n} \text{ is given for all }j]
 \le \left(\frac Cn\right)^{d-1}.$$
 Integrating in $t$ and interchanging the integral and the expectation, we obtain
$$
 \E[I_{2n}\,|E_j^{n} \text{ is given for all }j]\le \left(\frac Cn\right)^{d-1}.
$$
Therefore we can make choose squares from $n+1$ to step $2n$ such that $I_{2n} \leq (C/n)^{d-1}.$ We fix these squares and then the next steps are from $2n+1$ and $4n$. Starting from $n=1$, we estimated $I_n$ when $n$ is a power of $2$. Let $m\in \mathbb{N}$ be any step and choose a power of two $n$ such that $n\le m<2n$. We have 
\[
 P_t(E_j^m)\subset P_t(E_j^n).
\]
At step $n$ we have that $I_n\le(2C/n)^{d-1}$.
Multiplying over $j$, and integrating over $t$ gives
$I_m\le I_n$. Since $m<2n$, it follows that
\[
 I_m\le (2C/n)^{d-1}\le(4C/m)^{d-1}.
\]
Finally, for each sufficiently small $\delta$, choose $n$ with $4^{-n-1}<\delta\le4^{-n}$ and use  \eqref{slice} to get 
\[
 |K(\delta)|\lesssim_d I_n \lesssim (1/n)^{d-1}\lesssim \bigl(\log(1/\delta)\bigr)^{-(d-1)}.
\] This proves Theorem~\ref{thm:k}.
\paragraph{Argument for Probability one.}
Let's again fix the first $n$ many steps. For each chosen square \(Q_i\), enlarge \(Q_i\) to the unit square and let \(X_{Q_i}\) be the projected length of the set obtained after the next \(n\) steps. The argument \eqref{expectation} above, bounds $|P_t(E_j^{2n})|$ by the average of these $4^n$ independent lengths, i.e.
$$ |P_t(E_j^{2n})| \leq 4^{-n}\sum_{i=1}^{4^n}X_{Q_i} \qquad \text{ and } \qquad \E[4^{-n}\sum_{i=1}^{4^n}X_{Q_i}]\leq \frac Cn.$$
Let $A=C+1$. Using Hoeffding’s inequality, we get
\begin{align*}
   \Prob\{|P_t(E_j^{2n})|>A/n\mid E_j^{n} \text{ is given}\}
& \leq \Prob\{4^{-n}\sum_{i=1}^{4^n}X_{Q_i}>A/n\} \\& \leq \Prob\{4^{-n}\sum_{i=1}^{4^n}(X_{Q_i}-\E[X_{Q_i}])> 1/n\} \\
& \leq \exp\Big(\frac{-4^n}{2n^2}\Big). 
\end{align*}
We have used the trivial bound $X_{Q_i}\leq 2$ in the last inequality. Now let $B_n \in [0,1]$ be the set of $t$, where at least one of the projected lengths $|P_t(E_j^{2n})|$ is greater than $A/n$. We have
\begin{equation} \label{Integral}
  I_{2n} \leq \int_{t \in B_n}\prod_{j=1}^q |P_t(E_j^{2n})|\,dt+ \int_{t \neq B_n}\prod_{j=1}^q |P_t(E_j^{2n})|\,dt 2^{(d-1)} \leq \Big(\frac{A}{n}\Big)^{d-1}+2^{d-1}|B_n|.  
\end{equation}
We also have that 
$$ \Prob\{t \in B_n \mid E_j^{n} \text{ is given}\} \leq \sum_{i=1}^{d-1}\Prob\{|P_t(E_j^{2n})|>A/n\mid E_j^{n} \text{ is given}\} \leq (d-1)\exp\Big(\frac{-4^n}{2n^2}\Big)$$
Integrating  and taking expectation gives 
\begin{align*}|B_n| =\int_{0}^{1}1_{\{t \in B_n\}} \,dt\, \, \implies \E[|B_n|\mid E_j^{n} \text{ is given} ]& =\int_{0}^{1}\Prob\{t \in B_n \mid E_j^{n} \text{ is given}\} \, dt \\
&\leq (d-1)\exp\Big(\frac{-4^n}{2n^2}\Big).\end{align*}
Therefore, by Markov's inequality,
\[
 \Prob\{|B_n|>n^{-{(d-1)}}\mid E_j^{n} \text{ is given}\}\le  (d-1)n^{(d-1)}\exp\Big(\frac{-4^n}{2n^2}\Big).
\]
The right-hand side converges as $n$ goes over powers of 2. By the Borel-Cantelli lemma, almost surely $|B_n|\le n^{-(d-1)}$ for every sufficiently large such $n$. Therefore by \eqref{Integral}, $I_{2n}\le(A^{(d-1)}+2^{(d-1)})n^{-(d-1)}$, almost surely . The same nesting argument as above 
now gives the bound on $I_m$ for all sufficiently large $m \in \mathbb{N}$.
\end{proof}
\section{The truncated Nikodym set}\label{sec:nik}
Recall the definition \eqref{Trunc}, of truncated Nikodym set $N_{\epsilon}$, for $\epsilon >0$.
\begin{proof}[Proof of Theorem~\ref{thm:n}] We divide $[0,1]^2$ into thin strips of height $\epsilon$.
We will choose one horizontal line of points (and their corresponding line segments) in every other strip. We first introduce some notation for simplicity. \\
 At height $h$, the points belong to $N$ form the set
\[
 N_h=\{p\in[0,1]:(p,h)\in N\}.
\]
By Fubini's theorem, $|N_h|=0$ and $|[0,1]\setminus N_h|=1$ for almost every $h$. Let $D_h$ be the set 
\[
 D_h=\{(p,h):p\in[0,1]\setminus N_h\}.
\]
By definition of $N_{\epsilon}$, the line segment $(p,h)+t(z(p,h),1),\, \epsilon\leq t \leq 1$ lies in $N_{\epsilon}$. Let $\Phi_h$ be the map taking a point $(p,h)$ to its $x$-intercept and slope parameter of its line,
\[
 \Phi_h(p,h)=\bigl(p-hz(p,h),\,z(p,h)\bigr) \, \text{ and } E_N=\{\Phi_h(p,h) : (p,h) \in [0,1]^2\}.
\]
So $E_N$ is the set of line parameters (intercepts and slopes) through all the points in the
square. Let $(w,z)$ denote the $x$-intercept and slope parameter, its corresponding line meets height $y$ at $x$-coordinate $P_y(w,z)=w+yz$. \par 
Let $m=\lfloor1/(2\eps)\rfloor$, so $1/(4\eps)\le m\le1/(2\eps)$. For $0\le i<m$ pick $h_i\in(2i\eps,(2i+1)\eps)$ such that the line $\{p:(p,h_i)\in N\}$ has length zero $z_i(p)=z(p,h_i)$ is measurable. Since $|N_{\epsilon}|=0$, by Fubini's theorem almost all $h_i$ satisfy these two conditions.

Suppose $h_i<h_j$. A line coming from a point of $(p,h_i) \in D_{h_i}$ reaches height $h_j$ after moving up a distance $h_j-h_i$ from its point, i.e. $t=h_j-h_i$. We have chosen $h_i$ from alternate strips so for $i<j$ we have $\eps<h_j-h_i<1$. Thus
that point is in $N_{\epsilon}$ whereas $D_{h_j} \cap N =\varnothing.$ Therefore the two line families have no line in common, since a common line would meet height $h_j$ at the same point and we have 
\[
 \Phi_{h_i}(D_{h_i})\cap \Phi_{h_j}(D_{h_j})=\varnothing.
\]
At a fixed height $h_i$, the function $z_i(p)$ is measurable on the set $[0,1]\setminus N_{h_i}$ of length 1. By Lusin's theorem we can keep a compact subset $F_i\subset[0,1]\setminus N_{h_i}$ of length greater than $1/2$ on which $z_i$ is continuous. So $p-h_i z_i(p)$ is also continuous and we get the compact set \[
 P_i=\Phi_{h_i}(F_i\times\{h_i\})\subseteq\Phi_{h_i}(D_{h_i}).
\]
For two distinct heights $h_i$ and $h_j$, the corresponding compact sets $P_i$ and
$P_j$ are disjoint and thus the $L^1$ distance between them is non-zero. 
We choose $0<\kappa\le1$ small enough that
\[
 |w-w'|+|z-z'|\ge\kappa
 \qquad ((w,z)\in P_i,\ (w',z')\in P_j,\ i\ne j).
\]
Now fix $0<\delta<\min(1/16,\kappa/4)$. Take a maximal
$\delta$-separated set $p_{i,1}<\cdots<p_{i,n_i}$ in each $F_i$. Since $|F_i|\ge1/2$, we have that $n_i \sim \delta^{-1}$. Let $z_{i,k}=z_i(p_{i,k})$ and $w_{i,k}=p_{i,k}-h_i z_{i,k}$ and put
\[
 T_{i,k}=\{(x,y):h_i+\eps\le y\le h_i+1,\quad
                 |x-w_{i,k}-yz_{i,k}|<\delta\}.
\]
At each height $y \geq h_i+\eps$, the point $(w_{i,k}+yz_{i,k},y)$ belongs to the line segment in $N_\eps$. Every point of the tube segment $T_{i,k}$ is within distance $\delta$ of this point, so $T_{i,k}$ is in $N_\eps(\delta)$. Let the functions $f_i$ and $f$ be
\[
 f_i=\sum_k\mathbf1_{T_{i,k}},\qquad f=\sum_i f_i.
\]
Each tube $T_{i,k}$ has width $2\delta$ and height $1-\eps$. Hence
\begin{equation}\label{fi}
\int f_i=2\delta n_i(1-\eps)\gtrsim 1-\eps,
 \, \text{ and } \, \int f\gtrsim m(1-\eps).    
\end{equation}
Now we will control how much can the two strips overlap. Take two of the tubes $T$ and $T'$, with center lines $x=w+yz$ and $x=w'+yz'$ and assume they intersect at a height $y_0$. Then $|w-w'+y_0(z-z')| \leq 2\delta$. Thus, 
$$ \left|y-\frac{w'-w}{z-z'}\right|<\frac{2\delta}{|z-z'|}.$$
Thus, every point of \(T\cap T'\) lies in a range of heights of length at most \(4\delta/|z-z'|\). At each height the intersection has width at most \(2\delta\), so
\begin{equation}\label{eq:nik-strip-pair}
    |T\cap T'|\le \frac{8\delta^2}{|z-z'|}.
\end{equation}
Fix $h_i$ and consider the tubes through points $p_{i,k}$ and $p_{i,l}$. Let $r=|k-l|$. Suppose the Tubes meet at height $h_i+t$, where $\eps\le t\le1$. 
Thus
\[ \big|(p_{i,k}-p_{i,l})+t(z_{i,k}-z_{i,l})\big|<2\delta.\]
Using the triangle inequality, we have
\[
 r\delta\le |p_{i,k}-p_{i,l}|
 \le 2\delta+t|z_{i,k}-z_{i,l}|
 \le 2\delta+|z_{i,k}-z_{i,l}|.
\]
If $r\ge4$, then 
$|z_{i,k}-z_{i,l}| \geq r\delta/2$. Substituting this into \eqref{eq:nik-strip-pair} bounds the overlap of the tubes by $16\delta/r$.
For $r=1,2,3$, we use the trivial bound, the area of a whole tube, which is $2\delta(1-\eps)\le2\delta$. Hence $16\delta/r$ is an upper bound for every $r\ge1$. Therefore we have 
\begin{align}\label{combine}
\int f_i^2=\sum_k\sum_l |T_{i,k}\cap T_{i,l}|& \leq \sum_k \big(2\delta+16\delta\sum_{|l-k|\geq 1}^{}\frac1{|l-k|} \big)\notag \\ & \leq \big(2\delta+32\delta\sum_{r}^{n_i} \frac1{r} \big)
 \lesssim n_i\delta(1+\log n_i)\lesssim \log(1/\delta).    
\end{align}
Combining the above inequality with \eqref{fi} and then using Cauchy–Schwarz gives 
$$ 1\lesssim (\int f_i)^2\le|\bigcup_k T_{i,k}|\int f_i^2 \, \implies |\bigcup_k T_{i,k}| \geq (\log(1/\delta))^{-1}.$$
Now take two tubes through points $p_{i,k}$ and one based at $p_{j,k'}$, with $i\ne j$. Their line parameters belong to the disjoint sets $P_i$ and $P_j$. Let $\Delta w=w-w'$ and $\Delta z=z-z'$, we have $|\Delta w|+|\Delta z|\ge\kappa$ and $|\Delta w+y\Delta z|\leq 2\delta$. Suppose the tubes meet at height $y$. Since $h_i, t \in [0,1]$, so $0\le y\le2$. Therefore,
\[
\kappa \leq |\Delta w|+|\Delta z|=|(\Delta w+y\Delta z)-y\Delta z| +|\Delta z|
 \le |\Delta w+y\Delta z|+(y+1)|\Delta z|
 \le 2\delta+3|\Delta z|.
\]
This gives
$\kappa\le 2\delta+3|\Delta z|$. We have chosen $\delta<\kappa/4$, so $|\Delta z|\ge\kappa/6$.
By
\eqref{eq:nik-strip-pair} the intersection of the tubes has at most $48\delta^2/\kappa$. Therefore
\[
 \int f_i f_j
 =\sum_{k=1}^{n_i}\sum_{l=1}^{n_j}|T_{i,k}\cap T_{j,l}|
 \lesssim \frac{1}{\delta^2}\,\frac{\delta^2}{\kappa}
 =\frac{1}{\kappa}.
\] Combining this with \eqref{combine}, we get
\[
 \int f^2=\sum_i\int f_i^2+\sum_{i\neq j}\int f_i f_j
 \lesssim m(\log(1/\delta))+ m^2/\kappa.
\]
As $\delta$ decreases, $\log(1/\delta)$ increases while $m/\kappa$ stays fixed. For every
sufficiently small $\delta$ we therefore have $\log(1/\delta)\ge m/\kappa$, and the right hand side above is  $\lesssim m \log(1/\delta)$. The sufficiently small scale depends upon $\kappa$ and $m$ and thus upon $\epsilon$ and the set $N$. Applying
Cauchy–Schwarz and \eqref{fi} again, gives
$$ m^2\lesssim \Big(\int f\Big)^2 \leq |\bigcup_i\bigcup_k T_{i,k}| \int f^2 \lesssim m \log(1/\delta).$$ Therefore, 
$$ \epsilon^{-1}\log(1/\delta)^{-1} \lesssim |\bigcup_i\bigcup_k T_{i,k}| \leq N_\epsilon(\delta).$$
\end{proof}
\section{Kakeya families inside a Nikodym set}\label{sec:geo}
Fix a base height $h$ and consider the map
\[
 T_h(x,y)=\left(\frac{x}{y-h},\frac{1}{y-h}\right),\qquad y>h.
\]
On the segment through $(p,h)$ we have $y=h+t$ and $x=p+t z(p,h)$.
Dividing by $t$ gives $x/t=p/t+z(p,h)$. Thus, with $t'=1/t$ as the
new height, the image lies on the line $(pt'+z(p,h),t'),\, 1\leq t'\leq 2$.

Take $1/2\le t\le1$. Since $t=y-h$ is bounded below, $T_h$ is Lipschitz on each bounded part of this strip. The line segments lie in the null set $N$, so their images have zero measure. Also, $t'$ runs through $[1,2]$, giving a line in direction $(p,1)$ for each $p\in[0,1]$. These segments have length at least one and
fill an arc of directions. This give a planar Kakeya set of measure zero.

We can also see this connection in the overlap estimate of tubes. For a fixed height, if two tubes come from points \(r\) places apart in the \(\delta\)-separated list, their intersection has area at most \( \sim \delta/r\). Thus, when we add the overlaps within this family, we get the same sum \(\sum_{r\le 1/\delta}1/r\) and then the factor \(\log(1/\delta)\) as in the planar Kakeya argument.
For the bound in Theorem~\ref{thm:n} we use about $1/\eps$ of the horizontal line
families. There are about $1/\eps$ chosen lines. To obtain the
factor $1/\eps$, we needed to understand how much the families from different strips can overlap. What matters in the proof is that the chosen compact sets of line parameters are separated.
As $\delta$ tends to zero, the overlap between two different line families stays bounded. This is why all the chosen heights contribute to the lower bound.

There is also a density observation for the original punctured
segments, where $t$ can approach zero. Choose a countable dense set
of heights $H\subset[0,1]$. For each $h\in H$ and each
$p\in[0,1]$, the points $(p+t z(p,h),h+t)$ approach $(p,h)$ as $t\downarrow0$. Thus the closure of the union of these families contains $[0,1]\times H$, and hence the whole square. In this sense
the families are arranged densely along the square, and each gives a Kakeya family after the change of coordinates at its own height.

To see the difference of the line parameters from the Kakeya construction, consider the full punctured lines through the $[0,1]^2$.  Recall that $E_N=\{(p-hz(p,h),\,z(p,h) : (p,h) \in [0,1]^2\}$ and $P_t(w,z)=w+tz$. Thus
\[
 [0,1]\subseteq P_t(E_N)\qquad\text{for every }t\in[0,1].
\]
In the arc $0<\arctan t<\pi/4$ the projection of $E_n$ covers the whole
interval. In particular $\dim_{\rm H}E_N\ge1$. For \(t>1\), \(P_t(E_N)\) is contained in \(N\) at height \(t\). Since \(N\) has area zero, Fubini’s theorem gives \(|P_t(E_N)|=0\) for almost every \(t>1\).Thus on the second arc, $\pi/4<\arctan t<\pi/2$, the projection has length zero in almost every direction. By Marstrand's projection theorem $\dim_{\rm H}E_N=1$ and by Besicovitch Federer theorem we see that $E_N$ is purely unrectifiable and that $\mathcal H^1(E_N)=\infty$. Moreover $E_N$ is not $\sigma$-finite. \par
This gives a difference between the two planar constructions. In the Kakeya construction of Section~\ref{sec:cantor}, the line-parameter
set $E$ can be purely unrectifiable and have finite $\mathcal H^1$
measure. To obtain every slope $x\in[0,1]$, we only need the projection of $E$ onto the $x$-axis to contain $[0,1]$. For the Nikodym set, the points fill a square, so at \emph{every} height $t\in[0,1]$ the projection $P_t(E_N)$ contains
$[0,1]$. These heights give an interval of projection directions. It is this requirement that forces the purely unrectifiable set $E_N$ to have non-$\sigma$-finite $\mathcal H^1$ measure.

The projection argument also explains why Venetian blinds arise in the construction of Nikodym sets.  The two arcs explain the role of Venetian blinds. The Venetian-blind idea is to arrange many short pieces in a zigzag construction so that their projections still cover the required interval on the first arc, while increasing this overlap on the second, so that in the limit the projection has length zero for almost every direction on the second arc.

\section*{Acknowledgements}
I am thankful to Terence Tao for asking the two questions studied in this paper and for helpful discussions.

\textsc{Department of Mathematics, University of California, Los Angeles, 520 Portola Plaza, Los Angeles, CA 90095, USA.}
\textit{Email:}  \texttt{iqra@ucla.edu.}
\end{document}